\documentclass[a4paper,11pt]{article}
\usepackage[cp1251]{inputenc}
\usepackage{amssymb,amsmath,amsthm}
\author{Alexey\,V.\,Ustinov\footnote{The research of
the author was supported by Dynasty Foundation}}
\title{On the construction of a triangle from the feet of its angle
bisectors}

\begin{document}
\maketitle

Given a triangle $ABC$ with sides $a$, $b$, $c$. We want to construct a
triangle $A'B'C'$ such that segments $AA'$, $BB'$ and $CC'$ are its angle
bisectors. This problem was discussed in~\cite{Wernick1982,Meyers1996} (see
problem 138) and remained open. Article~\cite{Yiu2008} gives conic solution of
this problem.

Using formulae from~\cite{Yiu2008} it is easy to prove that in the general case
construction of a triangle from the feet of its angle bisectors with compass
and ruler is impossible. Following this article we denote by $(x:y:z)$
barycentric coordinates of incenter of triangle $A'B'C'$ with respect to
triangle $ABC$. Then vertices of triangle $A'B'C'$ will have coordinates
$(-x,y,z)$, $(x,-y,z)$, $(x,y,-z)$. As it was proved in~\cite{Yiu2008}, numbers
$x$, $y$, $z$ satisfy following equations:
\begin{align*}
-x(c^2y^2-b^2z^2)+yz((c^2+a^2-b^2)y-(a^2+b^2-c^2)z)=&0,\\
-y(a^2z^2-c^2x^2)+xz((a^2+b^2-c^2)z-(b^2+c^2-a^2)x)=&0,\\
-z(b^2x^2-a^2y^2)+xy((b^2+c^2-a^2)x-(c^2+a^2-b^2)y)=&0.\\
\end{align*}
Third equation follows from first and second because in our situation
$xyz\ne0$. Let $a=2$, $b=3$, $c=\sqrt{7}$ ($\angle C=60^{\circ}$). Using first
equation we can write $x$ in terms of $y$ and $z$. Substitutions into the
second equation lead to the irreducible equation for $t=z/(3y)$:
$$t^3+74t^2+259t-570=0.$$
It has three roots $t_1=1.52\ldots$, $t_2=-5.32\ldots$, $t_3=-70.19\ldots$.
First one gives the incenter of some triangle $A'B'C'$, another ones correspond
to excenters (of another triangles). It is impossible to construct $t_1$ using
compass and ruler (see \cite[ch.~3]{Kurant2001}), hence we can't construct the
incenter of triangle $A'B'C'$ and its vertices.

Author is greatful for V.\,Dubrovsky for the references
\cite{Meyers1996,Wernick1982}.

\end{document}